\documentclass[12pt,a4paper,draft]{article}
\usepackage{graphicx}
\usepackage{amsmath}
\usepackage{amsfonts}
\usepackage{amssymb}
\providecommand{\U}[1]{\protect\rule{.1in}{.1in}}
\newtheorem{theorem}{Theorem}

\newtheorem{definition}[theorem]{Definition}

\newtheorem{lemma}[theorem]{Lemma}

\newtheorem{proposition}[theorem]{Proposition}
\newtheorem{remark}[theorem]{Remark}

\newenvironment{proof}[1][Proof]{\noindent\textbf{#1.} }{\ \rule{0.5em}{0.5em}}

\begin{document}

\title{Infinite-dimensional generalization of Kolmogorov widths}
\author{O. Kounchev\\Institute of Mathematics and Informatics, \\Bulgarian Academy of Sciences\\and\\IZKS, University of Bonn\\Dedicated to the memory of Borislav Bojanov}
\maketitle
\begin{abstract}
Recently the theory of widths of Kolmogorov-Gelfand has received a great deal
of interest due to its close relationship with the newly born area of
Compressive Sensing in Signal Processing, cf. \cite{devore} and references
therein. However fundamental problems of the theory of widths in
multidimensional Theory of Functions remain untouched, as well as analogous
problems in the theory of multidimensional Signal Analysis. In the present
paper we provide a multidimensional generalization of the original result of
Kolmogorov about the widths of an "ellipsoidal sets" consisting of functions
defined on an interval.
\end{abstract}

\section{Introduction}

In his seminal paper \cite{kolmogorov1936} Kolmogorov has introduced the
theory of widths and applied it very successfully to the following set of
functions defined in the compact interval:
\begin{equation}
K_{p}:=\left\{  f\in AC^{p-1}\left(  \left[  a,b\right]  \right)  :%
{\displaystyle\int_{0}^{1}}
\left\vert f^{\left(  p\right)  }\left(  t\right)  \right\vert ^{2}%
dt\leq1\right\}  . \label{Kp}%
\end{equation}
In the present paper we consider a natural multivariate generalization of the
set $K_{p}$ given by
\begin{equation}
K_{p}^{\ast}:=\left\{  u\in H^{2p}\left(  B\right)  :%
{\displaystyle\int_{B}}
\left\vert \Delta^{p}u\left(  x\right)  \right\vert ^{2}dx\leq1\right\}  ,
\label{Kpstar}%
\end{equation}
where $\Delta^{p}$ is the $p-$th iterate of the Laplace operator $\Delta=%
{\displaystyle\sum_{j=1}^{n}}
\partial^{2}/\partial x_{j}^{2}$ in $\mathbb{R}^{n}.$ We generalize the notion
of width in the framework of the Polyharmonic Paradigm, and obtain analogs to
the one-dimensional results of Kolmogorov.

The \emph{Polyharmonic Paradigm} has been announced in \cite{okbook} as a new
approach in Multidimensional Mathematical Analysis (in particular, in the
Moment Problem, Approximation and Spline Theory) which is based on solutions
of higher order elliptic equations as opposed to the usual concept which is
based on algebraic and trigonometric polynomials of several variables. The
main result of the present research is a new aspect of the Polyharmonic
Paradigm. It provides a new hierarchy of infinite-dimensional spaces of
functions which are used for a generalization of the Kolmogorov's theory of
widths. This new hierarchy generalizes the hierarchy of finite-dimensional
subspaces $S_{N}$ of the space $C^{\infty}\left(  I\right)  $ for an interval
$I\subset\mathbb{R}.$ Let us give a \emph{rough idea} of this hierarchy in the
case of a domain $D\subset\mathbb{R}^{n},$ where $D$ is a compact domain with
sufficiently smooth boundary $\partial D.$ In the new hierarchy in
$\mathbb{R}^{n}$, the $N-$dimensional subspaces in $C^{\infty}\left(
I\right)  $ will be generalized by solution spaces
\[
S_{N}=\left\{  u:P_{2N}u\left(  x\right)  =0,\quad\text{for }x\in D\right\}
\subset C^{\infty}\left(  D\right)  ,
\]
where $P_{2N}$ is an elliptic operator of order $2N$ in the domain $D;$ their
precise definitions will be specified later on.

\section{Kolmogorov's result - a reminder}

Let us recall the original result of Kolmogorov provided in his seminal paper
\cite{kolmogorov1936} where he introduced for the first time the theory of
widths. Kolmogorov has considered the set $K_{p}$ defined in (\ref{Kp}). He
proved that this is an \textbf{ellipsoid} by constructing explicitly its
principal axes. Namely, he considered the eigenvalue problem
\begin{align}
\left(  -1\right)  ^{p}u^{\left(  2p\right)  }\left(  t\right)   &  =\lambda
u\left(  t\right)  \qquad\qquad\qquad\text{for }t\in\left(  0,1\right)
\label{eigen1}\\
u^{\left(  p+j\right)  }\left(  0\right)   &  =u^{\left(  p+j\right)  }\left(
1\right)  =0\qquad\text{for }j=0,1,...,p-1. \label{eigen2}%
\end{align}
By the results of M. Krein proved an year earlier \cite{krein}, \cite{naimark}%
, Kolmogorov proved that problem (\ref{eigen1})-(\ref{eigen2}) has the
following properties, cf. also \cite{lorentz}, Chapter $9.6,$ Theorem $9,$ p.
$146$, \cite{tikhomirov}, section $4.4.4,$ Theorem $6,$ p. $244$ ,
\cite{pinkus}, :

\begin{proposition}
Problem (\ref{eigen1})-(\ref{eigen2}) has a countable set of non-negative real
eigenvalues with finite multiplicity. If we denote them by $\lambda_{j}$ in a
monotone order, they satisfy $\lambda_{j}\longrightarrow\infty$ for
$j\longrightarrow\infty.$ They satisfy the following asymptotic $\lambda
_{j}=\pi^{2p}j^{2p}\left(  1+O\left(  j^{-1}\right)  \right)  .$ The
corresponding orthonormalized eigenfunctions $\left\{  \psi_{j}\right\}
_{j=1}^{\infty}$ form a complete orthonormal system in $L_{2}\left(  \left[
0,1\right]  \right)  .$ The eigenvalue $\lambda=0$ has multiplicity $p$ and
the corresponding eigenfunctions $\left\{  \psi_{j}\right\}  _{j=1}^{p}$ are
the basis for the solutions to equation $u^{\left(  p\right)  }\left(
t\right)  =0$ in the interval $\left(  0,1\right)  .$
\end{proposition}

Further, Kolmogorov provided a description of the axes of the "cylindrical
ellipsoid set" $K_{p}$, from which easily follows an approximation theorem of
\textbf{Jackson} \textbf{type}.

\begin{proposition}
\label{PKolmogorovJackson} Let $f\in L_{2}\left(  \left[  a,b\right]  \right)
$ have the $L_{2}-$expansion
\[
f\left(  t\right)  =%
{\displaystyle\sum_{j=1}^{\infty}}
f_{j}\psi_{j}\left(  t\right)  .
\]
Then $f\in K_{p}$ if and only if
\[%
{\displaystyle\sum_{j=1}^{\infty}}
f_{j}^{2}\lambda_{j}\leq1.
\]
For $N\geq p+1$ and every $f\in K_{p}$ holds the following estimate
(\textbf{Jackson type} approximation):
\[
\left\Vert f-%
{\displaystyle\sum_{j=1}^{N}}
f_{j}\psi_{j}\left(  t\right)  \right\Vert _{L_{2}}\leq\frac{1}{\sqrt
{\lambda_{N+1}}}=O\left(  \frac{1}{\left(  N+1\right)  ^{p}}\right)  .
\]
\end{proposition}

However, Komogorov didn't stop at this point but asked further, whether the
linear space $\widetilde{X}_{N}:=\left\{  \psi_{j}\right\}  _{j=1}^{N}$
provides the "best possible approximation among the linear spaces of dimension
$N$" in the following sense: if we put
\begin{equation}
d_{N}\left(  K_{p}\right)  :=\inf_{X_{N}}\operatorname*{dist}\left(
X_{N},K_{p}\right)  \label{KolmogorovWidth}%
\end{equation}
then Kolmogorov has proved in \cite{kolmogorov1936} the following equality%
\[
d_{N}\left(  K_{p}\right)  =\operatorname*{dist}\left(  \widetilde{X}%
_{N},K_{p}\right)  .
\]
Hence, the above result reads as
\begin{align*}
d_{N}\left(  K_{p}\right)   &  =\frac{1}{\sqrt{\lambda_{N+1}}}\qquad\text{for
}N\geq p\\
d_{N}\left(  K_{p}\right)   &  =\infty\qquad\qquad\text{for }N=0,1,...,p-1.
\end{align*}
Here we have used the notations
\begin{align*}
\operatorname*{dist}\left(  X,K_{p}\right)   &  :=\sup_{y\in K_{p}%
}\operatorname*{dist}\left(  X,y\right) \\
\operatorname*{dist}\left(  X,y\right)   &  =\inf_{x\in X}\left\Vert
x-y\right\Vert .
\end{align*}

\begin{definition}
The left quantity in (\ref{KolmogorovWidth}) is called \textbf{Kolmogorov
}$N-$\textbf{width}, while the best approximation space $\widetilde{X}_{N}$ is
called \textbf{extremal (optimal) subspace}, cf. \cite{lorentz},
\cite{tikhomirov}, \cite{pinkus}.
\end{definition}

Thus the \textbf{main concept of the theory of widths} is closely related to a
Jackson type theorem by which a special space $\widetilde{X}_{N}$ is
identified. Then one has to find in which sense is the space $\widetilde
{X}_{N}$ the extremal subspace. We may formulate it in other words: one has to
find as wide class of spaces $X_{N}$ as possible, among which $\widetilde
{X}_{N}$ is the extremal subspace.

Now let us consider the following set which is a \emph{natural multivariate
generalization} of the above set $K_{p}$ defined in (\ref{Kp}): For simplicity
sake we will restrict ourselves with the unit ball $\mathbb{B}$ in
$\mathbb{R}^{n}.$ We put
\[
K_{p}^{\ast}:=\left\{  u\in H^{2p}\left(  B\right)  :%
{\displaystyle\int_{B}}
\left|  \Delta^{p}u\left(  x\right)  \right|  ^{2}dx\leq1\right\}  .
\]
Let us remark that the \textbf{Sobolev space} $H^{2p}\left(  B\right)  $ is
the multivariate version of the space of absolutely continuous functions on
the interval with a highest derivative in $L_{2}$ (as in (\ref{Kp})). An
important feature of the set $K_{p}^{\ast}$ is that it contains an
infinite-dimensional subspace
\[
\left\{  u\in H^{2p}\left(  B\right)  :\Delta^{p}u\left(  x\right)
=0,\quad\text{for }x\in B\right\}  .
\]
Hence, all Kolmogorov widths are equal to infinity,
\[
d_{N}\left(  K_{p}^{\ast}\right)  =\infty\qquad\text{for }N\geq0
\]
and no way is seen to improve this if one remains within the
finite-dimensional setting.

The main purpose of the present paper is to find a proper setting in the
framework of the Polyharmonic Paradigm which generalizes the above results of Kolmogorov.

\section{Elliptic differential operators and Elliptic BVP}

As we said we restrict ourselves to a simple domain as the unit ball $B$ in
$\mathbb{R}^{n}.$ However the results below hold for a much bigger class of domains.

We will make extensive use of the following Green formula for the polyharmonic
operator $\Delta^{p},$ cf. \cite{aron}, p. $10$:\
\begin{equation}%
{\displaystyle\int_{B}}
\left(  \Delta^{p}u\cdot v-u\cdot\Delta^{p}v\right)  dx=%
{\displaystyle\sum_{j=0}^{p-1}}
{\displaystyle\int_{\partial B}}
\left(  \Delta^{j}u\cdot\partial_{n}\Delta^{p-1-j}v-\partial_{n}\Delta
^{j}u\cdot\Delta^{p-1-j}v\right)  ;\label{GreenPolyharmonic}%
\end{equation}
here $\partial_{n}$ denotes the normal derivative to $\partial B,$ for
functions $u$ and $v$ in the classes of Sobolev, $u,v$ $\in$ $H^{2p}\left(
B\right)  .$

For us the following eigenvalue problem will be important to consider for
$U\in H^{2p}\left(  B\right)  $:%
\begin{align}
\Delta^{2p}U\left(  x\right)   &  =\lambda U\left(  x\right)  \qquad\text{for
}x\in B\label{eigen1Multi}\\
\Delta^{p+j}U\left(  y\right)   &  =\partial_{n}\Delta^{p+j}U\left(  y\right)
=0,\qquad\text{for all }y\in\partial B,\quad j=0,1,...,p-1\label{eigen2Multi}%
\end{align}
where $\partial_{n}$ denotes the normal derivative at $y\in\partial B.$ The
operator $\Delta^{2p}$ is formally self-adjoint, cf. \cite{lions-magenes},
however the BVP (\ref{eigen1Multi})-(\ref{eigen2Multi}) is not a nice one from
the point of view of Elliptic Boundary Value problems. Since a direct
reference seems not to be available, we need a special consideration of this
problem provided in the following theorem.

\begin{theorem}
\label{TExpansionBerezanskii} Problem (\ref{eigen1Multi})-(\ref{eigen2Multi})
has only real non-negative eigenvalues.

1. The eigenvalue $\lambda=0$ has infinite multiplicity with corresponding
eigenfunctions $\left\{  \psi_{j}^{\prime}\right\}  _{j=1}^{\infty}$ which
represent an orthonormal basis of the space of all solutions to the equation
$\Delta^{p}U\left(  x\right)  =0,$ for $x\in B.$

2. The positive eigenvalues are countably many and each has \textbf{finite
multiplicity}, and if we denote them by $\lambda_{j}$ ordered increasingly,
they satisfy $\lambda_{j}\longrightarrow\infty$ for $j\longrightarrow\infty.$

3. The orthonormalized eigenfunctions, corresponding to eigenvalues
$\lambda_{j}>0,$ will be denoted by $\left\{  \psi_{j}\right\}  _{j=1}%
^{\infty}.$ The set of functions $\left\{  \psi_{j}\right\}  _{j=1}^{\infty
}\bigcup\left\{  \psi_{j}^{\prime}\right\}  _{j=1}^{\infty}$ form a complete
orthonormal system in $L_{2}\left(  B\right)  .$
\end{theorem}

\begin{remark}
Problem (\ref{eigen1Multi})-(\ref{eigen2Multi}) is widely known to be
non-regular elliptic BVP, as well as non-coercive variational, c.f.
\cite{agmon}, p. $150$ at the end of section $10$, Lions-Magenes Remark $9.8$
(chapter $2,$ section $9.6$, p. $240$ in the Russian edition) and section
$9.8$ there, p. $242.$ This problem will give us the eigenfunctions $\psi_{k}$
in the notations in \cite{lorentz}.
\end{remark}

The proof is provided in the Appendix below, section \ref{Sappendix}.

\section{The principal axes of the ellipsoid $K_{p}^{\ast}$ and Jackson type theorem}

Here we will find the principal exes of the ellipsoid $K_{p}^{\ast}$ defined
in (\ref{Kpstar}).

We prove the following theorem which generalizes Kolmogorov's one-dimensional
\cite{kolmogorov1936}, about the representation of the ellipsoid $K_{p}$ in
principal axes.

\begin{theorem}
Let $f\in K_{p}^{\ast}.$ Then $f$ is represented in a $L_{2}-$series as
\[
f\left(  x\right)  =%
{\displaystyle\sum_{j=1}^{\infty}}
f_{j}^{\prime}\psi_{j}^{\prime}\left(  x\right)  +%
{\displaystyle\sum_{j=1}^{\infty}}
f_{j}\psi_{j}\left(  x\right)
\]
(where by Theorem \ref{TExpansionBerezanskii} the eigenfunctions $\psi
_{j}^{\prime}$ satisfy $\Delta^{p}\psi_{j}^{\prime}\left(  x\right)  =0$ while
the eigenfunctions $\psi_{j}$ correspond to the eigenvalues $\lambda_{j}>0$ )
where the coefficients satisfy the inequality
\begin{equation}%
{\displaystyle\sum_{j=1}^{\infty}}
\lambda_{j}f_{j}^{2}\leq1.\label{EllipsoidCondition}%
\end{equation}
Vice versa, every sequence $\left\{  f_{j}^{\prime}\right\}  _{j=1}^{\infty
}\bigcup\left\{  f_{j}\right\}  _{j=1}^{\infty}$ with $%
{\displaystyle\sum_{j=1}^{\infty}}
\left|  f_{j}^{\prime}\right|  ^{2}+%
{\displaystyle\sum_{j=1}^{\infty}}
\left|  f_{j}\right|  ^{2}<\infty$ and $%
{\displaystyle\sum_{j=1}^{\infty}}
\lambda_{j}f_{j}^{2}\leq1$ define a function $f\in L_{2}\left(  B\right)  $
which is in $K_{p}^{\ast}.$
\end{theorem}%

\proof
1. According to Theorem \ref{TExpansionBerezanskii}, we know that arbitrary
$f\in L_{2}\left(  B\right)  $ is represented as
\begin{align*}
f\left(  x\right)   &  =%
{\displaystyle\sum_{j=1}^{\infty}}
f_{j}^{\prime}\psi_{j}^{\prime}\left(  x\right)  +%
{\displaystyle\sum_{j=1}^{\infty}}
f_{j}\psi_{j}\left(  x\right)  \\
\left\|  f\right\|  _{L_{2}}^{2} &  =%
{\displaystyle\sum_{j=1}^{\infty}}
\left|  f_{j}^{\prime}\right|  ^{2}+%
{\displaystyle\sum_{j=1}^{\infty}}
\left|  f_{j}\right|  ^{2}<\infty
\end{align*}
with convergence in the space $L_{2}\left(  B\right)  .$

2. From the proof of Theorem \ref{TExpansionBerezanskii}, we know that if we
put
\[
\phi_{j}\left(  x\right)  =\Delta^{p}\psi_{j}\left(  x\right)  \qquad\text{for
}j\geq1,
\]
then the system of functions
\[
\frac{\phi_{j}\left(  x\right)  }{\sqrt{\lambda_{j}}}\qquad\text{for }j\geq1
\]
is orthonormal sequence which is complete in $L_{2}\left(  B\right)  .$

3. We will prove now that if $f\in L_{2}\left(  B\right)  $ then $f\in
K_{p}^{\ast}$ iff
\[%
{\displaystyle\sum_{j=1}^{\infty}}
f_{j}^{2}\lambda_{j}\leq1.
\]
Indeed, for every $f\in H^{2p}\left(  B\right)  $ we have the expansion
$f\left(  x\right)  =%
{\displaystyle\sum_{j=1}^{\infty}}
f_{j}^{\prime}\psi_{j}^{\prime}\left(  x\right)  +%
{\displaystyle\sum_{j=1}^{\infty}}
f_{j}\psi_{j}\left(  x\right)  .$ We want to see that it is possible to
differentiate termwise this expansion, i.e.
\[
\Delta^{p}f\left(  x\right)  =%
{\displaystyle\sum_{j=1}^{\infty}}
f_{j}\Delta^{p}\psi_{j}\left(  x\right)  =%
{\displaystyle\sum_{j=1}^{\infty}}
f_{j}\phi_{j}\left(  x\right)
\]
Since $\left\{  \frac{\phi_{j}}{\sqrt{\lambda_{j}}}\right\}  _{j\geq1}$ is a
complete orthogonal basis of $L_{2}\left(  B\right)  $ it is sufficient to see
that
\[%
{\displaystyle\int_{B}}
\Delta^{p}f\left(  x\right)  \phi_{j}dx=%
{\displaystyle\int_{B}}
\left(
{\displaystyle\sum_{j=1}^{\infty}}
f_{j}\Delta^{p}\psi_{j}\left(  x\right)  \right)  \phi_{j}dx.
\]
Due to the boundary properties of $\phi_{j}$ and since $\phi_{j}=\Delta
^{p}\psi_{j},$ we obtain
\[%
{\displaystyle\int_{B}}
\Delta^{p}f\left(  x\right)  \phi_{j}dx=%
{\displaystyle\int_{B}}
f\left(  x\right)  \Delta^{p}\phi_{j}dx=\lambda_{j}%
{\displaystyle\int_{B}}
f\psi_{j}dx=\lambda_{j}f_{j}.
\]
On the other hand
\[%
{\displaystyle\int_{B}}
\left(
{\displaystyle\sum_{k=1}^{\infty}}
f_{k}\phi_{k}\left(  x\right)  \right)  \phi_{j}dx=\lambda_{j}f_{j}.
\]
Hence
\[
\Delta^{p}f\left(  x\right)  =%
{\displaystyle\sum_{j=1}^{\infty}}
f_{j}\Delta^{p}\psi_{j}\left(  x\right)  =%
{\displaystyle\sum_{j=1}^{\infty}}
f_{j}\phi_{j}\left(  x\right)  =%
{\displaystyle\sum_{j=1}^{\infty}}
\sqrt{\lambda_{j}}f_{j}\frac{\phi_{j}\left(  x\right)  }{\sqrt{\lambda_{j}}}%
\]
and since $\left\{  \frac{\phi_{j}}{\sqrt{\lambda_{j}}}\right\}  _{j\geq1}$ is
an orthonormal system, it follows
\[
\left\|  \Delta^{p}f\right\|  _{L_{2}}^{2}=%
{\displaystyle\sum_{j=1}^{\infty}}
\lambda_{j}f_{j}^{2}.
\]
Thus if $f\in K_{p}$ it follows that $%
{\displaystyle\sum_{j=1}^{\infty}}
\lambda_{j}f_{j}^{2}\leq1.$

Now, assume vice versa, that $%
{\displaystyle\sum_{j=1}^{\infty}}
f_{j}^{2}\lambda_{j}\leq1$ holds togather with $%
{\displaystyle\sum_{j=1}^{\infty}}
\left|  f_{j}^{\prime}\right|  ^{2}+%
{\displaystyle\sum_{j=1}^{\infty}}
\left|  f_{j}\right|  ^{2}<\infty$. We have to see that the function
\[
f\left(  x\right)  =%
{\displaystyle\sum_{j=1}^{\infty}}
f_{j}^{\prime}\psi_{j}^{\prime}\left(  x\right)  +%
{\displaystyle\sum_{j=1}^{\infty}}
f_{j}\psi_{j}\left(  x\right)
\]
belongs to the space $H^{2p}\left(  B\right)  .$ Based on the completeness and
orthonormality of the system $\left\{  \frac{\phi_{j}\left(  x\right)  }%
{\sqrt{\lambda_{j}}}\right\}  _{j=1}^{\infty}$ we may define the function
$g\in L_{2}$ by putting
\[
g\left(  x\right)  =%
{\displaystyle\sum_{j=1}^{\infty}}
\sqrt{\lambda_{j}}f_{j}\frac{\phi_{j}\left(  x\right)  }{\sqrt{\lambda_{j}}}=%
{\displaystyle\sum_{j=1}^{\infty}}
f_{j}\phi_{j}\left(  x\right)  ;
\]
it obviously satisfies $\left\|  g\right\|  _{L_{2}}\leq1.$

As is well known from the theory of Elliptic Boundary Value Problems we may
find a function $F\in H^{2p}\left(  B\right)  $ which is a solution to
equation $\Delta^{p}F=g$ (see Theorem $5.3$ in chapter $2,$ section $5.3,$
\cite{lions-magenes}). Let its representation be
\[
F\left(  x\right)  =%
{\displaystyle\sum_{j=1}^{\infty}}
f_{j}^{\prime}\psi_{j}^{\prime}\left(  x\right)  +%
{\displaystyle\sum_{j=1}^{\infty}}
F_{j}\psi_{j}\left(  x\right)
\]
with some $F_{j}$ satisfying $%
{\displaystyle\sum_{j}}
\left|  F_{j}\right|  ^{2}<\infty.$ As above we obtain
\begin{align*}
\lambda_{j}%
{\displaystyle\int_{B}}
F\psi_{j}dx &  =%
{\displaystyle\int_{B}}
F\Delta^{2p}\psi_{j}dx=%
{\displaystyle\int_{B}}
\Delta^{p}F\cdot\Delta^{p}\psi_{j}dx\\
&  =%
{\displaystyle\int_{B}}
g\cdot\phi_{j}dx
\end{align*}
which implies $F_{j}=f_{j}.$ Hence, $F=f$ and $f\in H^{2p}\left(  B\right)  .$
This ends the proof.%

\endproof

We are able to prove finally a \textbf{Jackson type} result as in Proposition
\ref{PKolmogorovJackson}.

\begin{theorem}
Let $N\geq1.$ Then for every $N\geq1$ and every $f\in K_{p}^{\ast}$ holds the
following estimate:
\[
\left\Vert f-%
{\displaystyle\sum_{j=1}^{\infty}}
f_{j}^{\prime}\psi_{j}^{\prime}\left(  x\right)  -%
{\displaystyle\sum_{j=1}^{N}}
f_{j}\psi_{j}\left(  x\right)  \right\Vert _{L_{2}}\leq\frac{1}{\sqrt
{\lambda_{N+1}}}.
\]
\end{theorem}%

\proof
The proof follows directly. Indeed, due to the monotonicity of $\lambda_{j},$
and inequality (\ref{EllipsoidCondition}), we obtain
\[
\left\Vert f-%
{\displaystyle\sum_{j=1}^{\infty}}
f_{j}^{\prime}\psi_{j}^{\prime}\left(  x\right)  -%
{\displaystyle\sum_{j=1}^{N}}
f_{j}\psi_{j}\left(  x\right)  \right\Vert _{L_{2}}^{2}=%
{\displaystyle\sum_{j=N+1}^{\infty}}
f_{j}^{2}\leq\frac{1}{\lambda_{N+1}}%
{\displaystyle\sum_{j=N+1}^{\infty}}
f_{j}^{2}\lambda_{j}\leq\frac{1}{\lambda_{N+1}}.
\]
This ends the proof.%

\endproof

Now we are able to prove a generalization of Kolmogorov's result about widths
\cite{kolmogorov1936}. It is important which classes of spaces we are going to
choose for generalizing the widths. We introduce the following subspaces in
$L_{2}\left(  B\right)  $: For integers $M\geq1$ we define
\begin{equation}
S_{M}:=\left\{  u\in H^{2M}\left(  B\right)  :Q_{2M}u\left(  x\right)
=0,\quad\text{for }x\in B\right\}  \label{SM}%
\end{equation}
where $Q_{2M}$ is a \emph{uniformly strongly elliptic} operator of order $2M,$
cf. \cite{agmondouglisNirenberg}, \cite{lions-magenes}, or \cite{okbook}, p.
$473$. We denote by $F_{N}$ a finite-dimensional subspace of $L_{2}\left(
B\right)  $ of dimension $N.$ We denote the special subspaces for
$P_{2M}=\Delta^{M}$ by
\begin{equation}
\widetilde{S}_{M}:=\left\{  u\in H^{2M}\left(  B\right)  :\Delta^{M}u\left(
x\right)  =0,\quad\text{for }x\in B\right\}  ,\label{StildeM}%
\end{equation}
and the special finite-dimensional subspaces
\begin{equation}
\widetilde{F}_{N}:=\left\{  \psi_{j}:j\leq N\right\}  _{lin}\label{FtildeN}%
\end{equation}
where $\psi_{j}$ are the eigenfunctions from Theorem
\ref{TExpansionBerezanskii}.

The following results are analogs to the original Kolmogorov's results about
widths, cf. \cite{kolmogorov1936}, or the more detailed exposition in
\cite{lorentz} (in Theorem $9,$ p. $146$), \cite{tikhomirov} and \cite{pinkus}.

\begin{theorem}
\label{TKolmogorovMultivariate} Let $Q_{2M}$ be a strongly elliptic
differential operator of order $2M$ in $B$, and let $N\geq0$ be arbitrary.

1. If $M<p$ then
\[
\operatorname*{dist}\left(  S_{M}%
{\textstyle\bigoplus}
F_{N},K_{p}^{\ast}\right)  =\infty.
\]
Hence,
\[
\inf_{Q_{2M}}\operatorname*{dist}\left(  S_{M}%
{\textstyle\bigoplus}
F_{N},K_{p}^{\ast}\right)  =\infty.
\]

2. If $M=p$ then
\[
\inf_{S_{p},F_{N}}\operatorname*{dist}\left(  S_{p}%
{\textstyle\bigoplus}
F_{N},K_{p}^{\ast}\right)  =\operatorname*{dist}\left(  \widetilde{S}_{p}%
{\textstyle\bigoplus}
\widetilde{F}_{N},K_{p}^{\ast}\right)  .
\]
\end{theorem}%

\proof
1. If we assume that $S_{M}$ and $\widetilde{S}_{p}$ are transversal the proof
is clear since $\widetilde{S}_{p}\subset K_{p}^{\ast}$ and there will be an
infinite-dimensional space in $\widetilde{S}_{p}\subset K_{p}^{\ast}$
containing infinite axes with direction $y\in\widetilde{S}_{p},$ such that
\[
\operatorname*{dist}\left(  S_{M}%
{\textstyle\bigoplus}
F_{N},y\right)  >0
\]
which implies
\[
\operatorname*{dist}\left(  S_{M}%
{\textstyle\bigoplus}
F_{N},K_{p}^{\ast}\right)  =\infty.
\]
If they are not transversal we apply Lemma \ref{LtransversalSpaces}; it is
clear that the finite-dimensional subspaces do not disturb the result, and the
proof is finished.

2. For proving the second item, let us first note that $\widetilde{S}%
_{p}\subset S_{p}%
{\textstyle\bigoplus}
F_{N}.$ Indeed, since $\widetilde{S}_{p}\subset K_{p}^{\ast}$ the violation of
$\widetilde{S}_{p}\subset S_{p}%
{\textstyle\bigoplus}
F_{N}$ would imply that there exists an infinite axis $y$ in $K_{p}^{\ast}$
not contained in $S_{p}%
{\textstyle\bigoplus}
F_{N}$ which would immediately give
\[
\operatorname*{dist}\left(  S_{p}%
{\textstyle\bigoplus}
F_{N},K_{p}^{\ast}\right)  =\infty.
\]
But by the Lemma \ref{LellipticCanonical} it follows that $P_{2p}=C\left(
x\right)  \Delta^{p}$ for some function $C\left(  x\right)  .$ Hence
$S_{p}=\widetilde{S}_{p}.$

Further we follow the usual way as in \cite{lorentz} to see that
$\widetilde{F}_{N}$ is extremal among all spaces $F_{N},$ i.e.
\[
\inf_{F_{N}}\operatorname*{dist}\left(  \widetilde{S}_{p}%
{\textstyle\bigoplus}
F_{N},K_{p}^{\ast}\right)  =\operatorname*{dist}\left(  \widetilde{S}_{p}%
{\textstyle\bigoplus}
\widetilde{F}_{N},K_{p}^{\ast}\right)  .
\]
This ends the proof.%

\endproof

We prove the following result which shows the mutual position of two subspaces:

\begin{lemma}
\label{LtransversalSpaces}Let the integers $M$ and $N$ satisfy $M<N,$ and the
integer $M_{1}\geq0.$ Then for the corresponding $S_{M}$ and $S_{N}$ defined
in (\ref{SM}) by the operators $P_{2M}$ and $Q_{2N}=\Delta^{N}$ respectively,
holds
\[
\operatorname*{dist}\left(  S_{M}%
{\textstyle\bigoplus}
F_{M_{1}},S_{N}\right)  =\infty.
\]
There is a linear subspace $Y_{N-M}\subset S_{N}$ with $Y_{N-M}\perp S_{M}$
and it is an infinite-dimensional space of solutions to an Elliptic Boundary
Value Problem.
\end{lemma}%

\proof
Let us consider the case $M_{1}=0.$ For the uniformly strongly elliptic
operator $P_{2M}$ we choose the Dirichlet system of boundary operators
$B_{j}=\frac{\partial^{j-1}}{\partial n^{j-1}}.$ It is a classical fact (cf.
Remark $1.3$ in chapter $2,$ section $1.4$ in \cite{lions-magenes}) that this
system satisfies conditions (iii) in section $5.1,$ chapter $2$ in
\cite{lions-magenes}, or in other words, the system of operators $\left\{
P_{2M};\frac{\partial^{j}}{\partial n^{j}}:j=0,1,...,M-1\right\}  $ forms a
\emph{regular Elliptic Boundary Value Problem }(thisi is the so-called called
Dirichlet BVP associated with the operator $P_{2M}$). \emph{ }Hence, we may
apply the existence Theorem $5.2$ and Theorem $5.3$ in \cite{lions-magenes}.
As in Theorem $2.1$ (section $2.2,$ chapter $2$ in \cite{lions-magenes}) we
complete the system $\left\{  B_{j}\right\}  _{j=1}^{M}$ by the system of
boundary operators $S_{j}=\frac{\partial^{M-1+j}}{\partial n^{M-1+j}}.$ Hence,
the system composed $\left\{  B_{j}\right\}  _{j=1}^{M}\bigcup\left\{
S_{j}\right\}  _{j=1}^{M}$ is a Dirichlet system of order $2M$ (cf. e.g.
Definition $23.12,$ p. $474$ in \cite{okbook}). Further, by Theorem $2.1$ in
\cite{lions-magenes} quoted above, there exists a unique Dirichlet system of
order $2M$ of boundary operators $\left\{  C_{j},T_{j}\right\}  _{j=1}^{M}$
which is uniquely determined as the adjoint to the system $\left\{
B_{j},S_{j}\right\}  _{j=1}^{M},$ and the following Green formula holds:
\begin{equation}%
{\displaystyle\int_{B}}
\left(  P_{2M}u\cdot v-u\cdot P_{2M}^{\ast}v\right)  dx=%
{\displaystyle\sum_{j=1}^{M}}
{\displaystyle\int_{\partial B}}
\left(  S_{j}u\cdot C_{j}v-B_{j}u\cdot T_{j}v\right)  d\sigma_{y}%
,\label{GreenGeneral}%
\end{equation}
for all $u,v\in H^{2M}\left(  B\right)  ;$ here $d\sigma_{y}$ denotes the
surface element on the sphere $\partial B.$

We consider the elliptic operator $\Delta^{N}P_{2M}.$ As a product of two
\emph{uniformly strongly elliptic} operators it is such again. By a standard
construction of Theorem $2.1$ in \cite{lions-magenes} cited above (section
$2.2,$ chapter $2$ in \cite{lions-magenes}), we complete the Dirichlet system
of operators $\left\{  B_{j},S_{j}\right\}  _{j=1}^{M}$ with $N-M$ boundary
operators $R_{j}=\frac{\partial^{2M-1+j}}{\partial n^{2M-1+j}},$
$j=1,2,...,N-M.$ Again by the above cited theorem, the Dirichlet system of
boundary operators
\[
\left\{  B_{j},S_{j}\right\}  _{j=1}^{M}%
{\textstyle\bigcup}
\left\{  R_{j}\right\}  _{j=1}^{N-M}%
\]
covers the operator $\Delta^{N}P_{2M}.$ Finally, we consider the solutions
$g\in H^{2N+2M}\left(  B\right)  $ to the following Elliptic Boundary Value
Problem:
\begin{align}
\Delta^{N}P_{2M}^{\ast}g\left(  x\right)   &  =0\qquad\qquad\qquad
\ \ \text{for }x\in B\label{gorthogonal1}\\
B_{j}g\left(  y\right)   &  =S_{j}g\left(  y\right)  =0\qquad\text{for
}j=0,1,...,N-1,\text{ for }y\in\partial B\label{gorthogonal2}\\
R_{j}g\left(  y\right)   &  =h_{j}\left(  y\right)  \qquad\qquad\ \ \text{for
}j=1,2,...,N-M,\text{ for }y\in\partial B.\label{gorthogonal3}%
\end{align}
We may apply the existence Theorem $5.2$ and Theorem $5.3$ in chapter $2$ in
\cite{lions-magenes}, to the solvability of problem (\ref{gorthogonal1}%
)-(\ref{gorthogonal3}) in the space $H^{2M+2N}\left(  B\right)  .$

First of all, it is clear from (\ref{gorthogonal1}) that $P_{2M}^{\ast}g\in
S_{N}.$

Let us check the properties of the function $P_{2M}^{\ast}g.$ By the Green
formula (\ref{GreenGeneral}), the function $P_{2M}^{\ast}g$ satisfies
$P_{2M}^{\ast}g\perp S_{M}$ , or equivalently,
\[%
{\displaystyle\int_{B}}
P_{2M}^{\ast}g\cdot vdx=0\qquad\text{for all }v\text{ with }P_{2M}v=0.
\]

By the general existence Theorem $5.3$ (the Fredholmness property) in
\cite{lions-magenes} mentioned above, we know that a solution $g$ to problem
(\ref{gorthogonal1})-(\ref{gorthogonal3}) exists for those boundary data
$\left\{  h_{j}\right\}  _{j=1}^{N-M}$ which satisfy only a finite number of
linear restrictions, provided by conditions (5.18) there; these are determined
by the solutions to the homogeneous adjoint Elliptic BVP. Hence, it follows
that the set $Y_{N-M}$ of the functions $P_{2M}^{\ast}g$ where $g$ is a
solution to (\ref{gorthogonal1})-(\ref{gorthogonal3}) is infinite-dimensional.
It follows that the space $S_{N}\setminus S_{M}$ is infinite-dimensional as
well, hence
\[
\operatorname*{dist}\left(  S_{M}%
{\textstyle\bigoplus}
F_{M_{1}},S_{N}\right)  =\infty.
\]

Since obviously a fininte-dimensional subspace $F_{M_{1}}$ would not disturb
the above argumentation, this ends the proof.%

\endproof

\begin{remark}
Lemma \ref{LtransversalSpaces} may be considered as a generalization in our
setting of a theorem of Gohberg-Krein of $1957$ (cf. \cite{lorentz}, Theorem
$2$ on p. $137$ ) in a Hilbert space.
\end{remark}

We need the following intuitive result which is however not trivial.

\begin{lemma}
\label{LellipticCanonical}Let for some elliptic differential operator $P_{2N}$
of order $2N$ the following inclusion hold $S_{N}\subset\widetilde{S}%
_{N}\setminus F,$ i.e.
\begin{align*}
&  \left\{  u\in H^{2N}\left(  B\right)  :P_{2N}u\left(  x\right)  =0,\quad
x\in B\right\}  \subset\\
&  \subset\left\{  u\in H^{2N}\left(  B\right)  :\Delta^{N}u\left(  x\right)
=0,\quad x\in B\right\}  \setminus F,
\end{align*}
where $F\subset L_{2}\left(  B\right)  $ is a finite-dimensional space. Then
\begin{equation}
P_{2N}\left(  x,D_{x}\right)  =c\left(  x\right)  \Delta^{N}\label{P2N=cx}%
\end{equation}
for some function $C\left(  x\right)  .$
\end{lemma}%

\proof
Since the general case is rather technical we will consider only $N=1$ in
$B\subset\mathbb{R}^{2}.$ It is clear that the arguments are purely local so
we will prove that equality (\ref{P2N=cx}) holds at $\left(  x_{1}%
,x_{2}\right)  =x=0\in B.$ Assume that
\[
P_{2N}\left(  x,D_{x}\right)  u\left(  x\right)  =a\left(  x\right)
u_{x_{1},x_{1}}+2b\left(  x\right)  u_{x_{1},x_{2}}+c\left(  x\right)
u_{x_{2},x_{2}}+d\left(  x\right)  u_{x_{1}}+e\left(  x\right)  u_{x_{2}%
}+f\left(  x\right)  u;
\]
here $w_{x_{j}}$ denotes the partial derivative $\frac{\partial w}{\partial
x_{j}}.$ By assumption, for the function $u\in\widetilde{S}_{1}\setminus F$
holds also
\[
\left(  a\left(  x\right)  -c\left(  x\right)  \right)  u_{x_{1},x_{1}%
}+2b\left(  x\right)  u_{x_{1},x_{2}}+d\left(  x\right)  u_{x_{1}}+e\left(
x\right)  u_{x_{2}}+f\left(  x\right)  u=0.
\]
Let us denote the following harmonic functions by $u^{j}$ for $j=1,2,...,6,$
as follows: $1,$ $x_{1},$ $x_{2},$ $x_{1}^{2}-x_{2}^{2},$ $x_{1}x_{2}$ . Let
us assume that they do not belong to $F.$ We see that the Jacobi matrix of
these functions at $x_{1}=x_{2}=0,$ is
\[
\left(
\begin{array}
[c]{ccccc}%
u_{x_{1},x_{1}}^{j} & u_{x_{1},x_{2}}^{j} & u_{x_{1}}^{j} & u_{x_{2}}^{j} &
u^{j}%
\end{array}
\right)  _{j=1}^{5}=\left(
\begin{array}
[c]{ccccc}%
0 & 0 & 0 & 0 & 1\\
0 & 0 & 1 & 0 & 0\\
0 & 0 & 0 & 1 & 0\\
2 & 0 & 0 & 0 & 0\\
0 & 1 & 0 & 0 & 0
\end{array}
\right)
\]
which is obviously non-degenerate. Hence, $a\left(  0\right)  -c\left(
0\right)  =b\left(  0\right)  =d\left(  0\right)  =e\left(  0\right)
=f\left(  0\right)  =0.$

In the case if some of the above functions $u^{j}$ belongs to the space $F,$
it is possible to approximate it by other harmonic functions also including up
to their second derivatives at $0$ (one may apply approximation arguments as
in \cite{hedberg}). The respective Jacobian will be non-zero and the
conclusion of the theorem will follow. This ends the proof.%

\endproof

The proof of Theorem \ref{TKolmogorovMultivariate} above permits a much bigger
generalization which will be provided in a forthcoming paper.

\section{Appendix on Elliptic Boundary Value Problems \label{Sappendix}}

\subsection{Proof of Theorem \ref{TExpansionBerezanskii}}%

\proof
\textbf{(1)} We consider the following auxiliary elliptic \emph{eigenvalue
problem}
\begin{align}
\Delta^{2p}\phi\left(  x\right)   &  =\lambda\phi\left(  x\right)
\qquad\qquad\text{on }B,\label{ffiBVP1}\\
\partial\Delta^{j}\phi\left(  y\right)   &  =\Delta^{j}\phi\left(  y\right)
=0\qquad\text{for }j=0,1,...,p-1,\text{ for }y\in\partial B.\label{ffiBVP2}%
\end{align}
It is straigthforward to check that this is a \emph{regular Elliptic BVP}
considered in the Sobolev space $H^{2p}\left(  B\right)  $ since it satisfies
all conditions (i)-(iii) in chapter $2,$ section $5.1$, \cite{lions-magenes},
cf. also \cite{egorov-shubin0}. Hence, we are able to apply the existence
theorems in section $5.3$ there. Further, it is straightforward to check that
it is a self-adjoint problem (cf. section $2.5,$ chapter $2,$
\cite{lions-magenes}): in the polyharmonic Green formula
(\ref{GreenPolyharmonic}) we put $\left\{  B_{j}\right\}  _{j=1}^{2p}=\left\{
\partial\Delta^{j},\Delta^{j}\right\}  _{j=0}^{p-1}$ and we see that in the
context of the general Green formula (\ref{GreenGeneral}) the adjoint system
of operators $\left\{  C_{j}\right\}  _{j=1}^{2p}=\left\{  \partial\Delta
^{j},\Delta^{j}\right\}  _{j=0}^{p-1}$ which proves the self-adjointness of
problem (\ref{ffiBVP1})-(\ref{ffiBVP2}). Hence, we may apply the main results
about the Spectral theory of regular self-adjoint Elliptic BVP. We refer to
\cite{egorov-shubin0}, section $3$ in chapter $2,$ p. $122,$ Theorem $2.52,$
and to references therein (cf. in particular the monograph of Yu. Berezanskii
devoted to expansions in eigenfunctions \cite{berezanskii}, chapter $6,$
section $2$).

By the uniqueness Lemma \ref{LuniquenessDirichlet} the eigenvalue problem
(\ref{ffiBVP1})-(\ref{ffiBVP2}) has only zero solution for $\lambda=0.$ It has
eigenfunctions $\phi_{k}\in H^{2p}\left(  B\right)  $ with eigenvalues
$\lambda_{k}>0$ for $k=1,2,3,...$ for which $\lambda_{k}\longrightarrow\infty$
as $k\longrightarrow\infty.$

\textbf{(2)} Next we consider the problem
\begin{align}
\Delta^{2p}\varphi &  =\phi_{k}\label{fiBVP1}\\
\partial\Delta^{j}\varphi\left(  y\right)   &  =\Delta^{j}\varphi\left(
y\right)  =0\qquad\text{for }j=0,1,...,p-1,\text{ for }y\in\partial
B,\label{fiBVP2}%
\end{align}
in the Sobolev space $H^{2p}\left(  B\right)  .$ Obviously, the Elliptic BVP
defined by problem (\ref{fiBVP1})-(\ref{fiBVP2}) coincides with the Elliptic
BVP defined by (\ref{ffiBVP1})-(\ref{ffiBVP2}) and all remarks there hold.
Hence, problem (\ref{fiBVP1})-(\ref{fiBVP2}) has \textbf{unique solution}
$\varphi_{k}\in H^{2p}\left(  B\right)  .$ We put
\[
\psi_{k}=\Delta^{p}\varphi_{k}.
\]
Hence, $\Delta^{p}\psi_{k}=\phi_{k}.$ We infer that on the boundary $\partial
B$ hold the equalities $\Delta^{p+j}\psi_{k}=\Delta^{j}\phi_{k}$ and
$\partial\Delta^{p+j}\psi_{k}=\partial\Delta^{j}\phi_{k};$ since $\phi_{k}$
are solutions to (\ref{ffiBVP1})-(\ref{ffiBVP2}) it follows
\begin{equation}
\Delta^{p+j}\psi_{k}\left(  y\right)  =\partial\Delta^{p+j}\psi_{k}\left(
y\right)  =0\qquad\text{for }j=0,1,...,p-1,\text{ for }y\in\partial
B.\label{psiBVP}%
\end{equation}

We will prove that $\psi_{k}$ are solutions to problem (\ref{eigen1Multi}%
)-(\ref{eigen2Multi}), they are mutulally \textbf{orthogonal,} and they are
also orthogonal to the space $\left\{  v\in H^{2p}:\Delta^{p}v=0\right\}  $.

Let us see that
\[
\Delta^{2p}\psi_{k}=\lambda_{k}\psi_{k}.
\]
By the definition of $\psi_{k}$ this is equivalent to
\[
\Delta^{3p}\varphi_{k}=\lambda_{k}\Delta^{p}\varphi_{k};
\]
from $\Delta^{2p}\varphi_{k}=\phi_{k}$ this is equivalent to
\[
\Delta^{p}\phi_{k}=\lambda_{k}\Delta^{p}\varphi_{k}%
\]
On the other hand, we have obviously $\Delta^{2p}\phi_{k}=\lambda_{k}%
\Delta^{2p}\varphi_{k}$ by the basic properties of $\phi_{k}$ and $\varphi
_{k},$ hence
\[
\Delta^{2p}\left(  \phi_{k}-\lambda_{k}\varphi_{k}\right)  =0.
\]
Note that both $\phi_{k}$ and $\varphi_{k}$ sastisfy the same zero boundary
conditions, namely (\ref{ffiBVP2}) and (\ref{fiBVP2}). Hence, by the
uniqueness Lemma \ref{LuniquenessDirichlet} it follows that $\phi_{k}%
-\lambda_{k}\varphi_{k}=0$ which implies $\Delta^{2p}\psi_{k}=\lambda_{k}%
\psi_{k}.$ Thus we see that $\psi_{k}$ is a solution to problem
(\ref{eigen1Multi})-(\ref{eigen2Multi}) and does not satisfy $\Delta^{p}%
\psi=0$ !

The orthogonality to the subspace $\left\{  v\in H^{2p}:\Delta^{p}v=0\right\}
$ follows easily from the Green formula (\ref{GreenPolyharmonic}) and the zero
boundary conditions (\ref{psiBVP}) of $\psi_{k},$ by the following:
\begin{align*}
&
{\displaystyle\int_{D}}
\left(  \Delta^{2p}\psi_{k}\cdot v-\psi_{k}\cdot\Delta^{2p}v\right)  dx\\
&  =%
{\displaystyle\sum_{j=0}^{2p-1}}
{\displaystyle\int_{\partial D}}
\left(  \Delta^{j}\psi_{k}\cdot\partial_{n}\Delta^{2p-1-j}v-\partial_{n}%
\Delta^{j}\psi_{k}\cdot\Delta^{2p-1-j}v\right)
\end{align*}
and since $%
{\displaystyle\int_{D}}
\Delta^{2p}\psi_{k}\cdot vdx=\lambda_{k}%
{\displaystyle\int_{D}}
\psi_{k}\cdot vdx.$

The orthonormality of the system $\left\{  \psi_{k}\right\}  _{k=1}^{\infty}$
follows now easily by the equality
\[
\lambda_{k}%
{\displaystyle\int}
\psi_{k}\psi_{j}dx=%
{\displaystyle\int}
\Delta^{2p}\psi_{k}\psi_{j}dx=%
{\displaystyle\int}
\Delta^{p}\psi_{k}\Delta^{p}\psi_{j}dx=%
{\displaystyle\int}
\phi_{k}\phi_{j}dx
\]
and the orthogonality of the system $\left\{  \phi_{k}\right\}  _{k=1}%
^{\infty}.$ For the completeness of the system $\left\{  \psi_{k}\right\}
_{k=1}^{\infty}$, let us assume that for some $f\in L_{2}\left(  B\right)  $
holds
\begin{equation}%
{\displaystyle\int_{B}}
f\cdot\psi_{k}dx=%
{\displaystyle\int_{B}}
f\cdot\psi_{k}^{\prime}dx=0\qquad\text{for all }k\geq1.\label{fOrthogonal}%
\end{equation}
Then the Green formula (\ref{GreenPolyharmonic}) implies
\begin{align*}
0 &  =\lambda_{k}%
{\displaystyle\int_{B}}
f\cdot\psi_{k}dx=%
{\displaystyle\int_{B}}
f\cdot\Delta^{2p}\psi_{k}dx=%
{\displaystyle\int_{B}}
\Delta^{p}f\cdot\Delta^{p}\psi_{k}dx\\
&  =%
{\displaystyle\int_{B}}
\Delta^{p}f\cdot\phi_{k}dx\qquad\text{for all }k\geq1.
\end{align*}
By the completeness of the system $\left\{  \phi_{k}\right\}  _{k\geq1}$ this
implies that $\Delta^{p}f=0.$ From the second orthogonality in
(\ref{fOrthogonal}) follows that $f\equiv0,$ and this ends the proof of the
completeness of the system $\left\{  \psi_{j}^{\prime}\right\}  _{j=1}%
^{\infty}\bigcup\left\{  \psi_{j}\right\}  _{j=1}^{\infty}.$%

\endproof

We have used above the following simple result.

\begin{lemma}
\label{LuniquenessDirichlet} The solution to problem (\ref{ffiBVP1}%
)-(\ref{ffiBVP2}) for $\lambda=0$ is unique.
\end{lemma}%

\proof
From Green formula (\ref{GreenPolyharmonic}) we obtain
\[%
{\displaystyle\int_{B}}
\left[  \Delta^{p}\phi\right]  ^{2}dx=%
{\displaystyle\int}
\phi\cdot\Delta^{2p}\phi dx=0,
\]
hence $\Delta^{p}\phi=0.$ Now we apply the second Green formula (2.11) in
\cite{aron} which infers immediately $\phi\equiv0.$%

\endproof

\textbf{Acknowledgement:} The author acknowldges the support of the Alexander
von Humboldt Foundation, and of Project Astroinformatics, DO-02-275 with
Bulgarian NSF. The author thanks Prof. Matthias Lesch for the interesting
discussion about hierarchies of infinite-dimensional linear spaces. I have got
a good advice on the elliptic BVP (\ref{eigen1Multi})-(\ref{eigen2Multi}) from
a conversation with Prof. P. Popivanov, N. Kutev and D. Boyadzhiev.

\end{document}